\documentclass[12pt,a4]{article}
\usepackage{amsfonts}
\usepackage{amssymb}
\def\fdem {{$\hfill\square$}}
\title{  Hartogs-Bochner type theorem in projective space}
\author{Fr\'ed\'eric Sarkis$^*$}
\date{ }
\newtheorem{defi}{Definition}
\newtheorem{theo}{Theorem}
\newtheorem{prop}{Proposition}
\newtheorem{lemme}{Lemma}
\newtheorem{coro}{Corollary}

\newenvironment{demo}{{\em Proof.}}{}
\newtheorem{problem}{Problem}
\begin{document}
\maketitle

\begin{center}
\begin{minipage}[t]{13.5cm}
{\bf abstract.}
{\footnotesize 
We prove the following Hartogs-Bochner type theorem:
{\it Let $M$ be a connected $C^2$ 
hypersurface of $P_n(\mathbb{C})$ ($n\geq 2$) which divides 
$P_n(\mathbb{C})$ in two connected open sets $\Omega_1$ and $\Omega_2$.
Then there exists $i \in \{1,2\}$ such that
$C^{1}$ CR 
functions defined on $M$ extends holomorphically to $\Omega_i$.}}
\end{minipage}
\end{center}

{\let\thefootnote\relax
\footnote{\hskip -2mm{\it *} : Supported by the TMR network} 

\section{Introduction}

The classical Hartogs-Bochner theorem states that 
if $\Omega \subset \subset \mathbb{C}^n$ ($n\geq 2$) is
a domain which boundary $\partial  \Omega$ is smooth
and connected, then every continuous CR function defined
on $\partial \Omega$ extends holomorphically to $\Omega$.
A natural question is to ask if such an extension phenomenon
is valid for domains included in a complex manifold $X$.
Of course, in the case that $X$ is compact,
there is no hope to expect such a result. Indeed,
if the Hartogs-Bochner phenomenon is valid in $X$,
then CR functions on $\partial \Omega$ would extend
to $\Omega$ but also to $X \backslash \Omega$ and thus
are constant which is impossible in general.
Nevertheless, the following Hartogs-Bochner type phenomenon
has been conjectured in $P_2(\mathbb{C})$:
{\it Let $M$ be a connected $C^2$ hypersurface
of  $P_2(\mathbb{C})$ which divides $P_2(\mathbb{C})$
into two connected open sets $\Omega_1$ and $\Omega_2$. 
Then CR functions on $M$ extend holomorphically 
to one of these sets.}\par
This conjecture has interested many authors at least since 1996
when E. Porten communicated to me the question with reference to
R. Dwilewicz. In \cite{sarkis}, we proved that holomorphic
(resp. meromorphic) functions defined in a connected
neighborhood of $M$ extend holomorphically (resp. meromorphically)
to one of the two sides of $M$ and repeated the question about
the extension of CR functions. Recently, Dwilewicz and Merker \cite{dwil}
gave a simplification of this prove in the holomorphic case and
raised again the question.
In \cite{iordan}, Henkin and Iordan gave a proof
of the conjecture for $M$ of Lipschitz class
but only under the hypothesis that one of the
two sides of $M$ contains a weakly concave domain with smooth
boundary. \par

In this paper, we prove the following Hartogs-Bochner type theorem:
\begin{theo}
Let $M$ be a connected $C^2$ hypersurface of $P_n(\mathbb{C})$ 
($n\geq 2$) which divides 
$P_n(\mathbb{C})$ in two connected open sets $\Omega_1$ and $\Omega_2$.
Then there exists $i \in \{1,2\}$ such that
for any $C^1$ CR function $f:M \rightarrow \mathbb{C}$,
there exists a holomorphic function $F \in {\mathcal{O}}(\Omega_i)
\cap C^1(\overline \Omega_i)$ such that $F|_M=f$.
\end{theo}
Moreover, in the case that there exists a non-constant
$C^1$ CR function $g$ on $M$,
we can give a more precise statement:
\begin{theo}
Let $M$ be a connected $C^2$ hypersurface of $P_n(\mathbb{C})$ 
($n\geq 2$) which divides 
$P_n(\mathbb{C})$ in two connected open sets $\Omega_1$ and $\Omega_2$.
Let us suppose that there exists a non-constant $C^1$ CR function
$g: M \rightarrow \mathbb{C}$, then there exists 
$i \in \{1,2\}$ such that:
\begin{enumerate}
\item The CR function $g$ admits a holomorphic ($C^1$ up to the boundary)
extension to $\Omega_i$.
\item For every continuous CR function $f:M \rightarrow \mathbb{C}$,
there exists a function $F \in \mathcal{O}(\Omega_i) 
\cap C^0(\overline \Omega_i)$ such that $F|_M=f$.
\item ${\mathcal{O}}(\Omega_{j}) \cap C^0(\overline \Omega_{j}) 
\equiv \mathbb{C}$ where $\Omega_j=P_n(\mathbb{C}) \backslash
\overline \Omega_i$.
\end{enumerate}
\end{theo}

If we suppose that $n=2$, then theorem $2$
is also valid for smooth CR maps:
\begin{theo}
Let $M$ be a connected $C^2$ hypersurface of $P_2(\mathbb{C})$ 
which divides 
$P_2(\mathbb{C})$ in two connected open sets $\Omega_1$ and $\Omega_2$.
Let us suppose that there exists a non-constant $C^1$ CR function
$g: M \rightarrow \mathbb{C}$, then there exists
 $i \in \{1,2\}$ such that
every $C^1$ CR map $f: M \rightarrow Y$ (where $Y$
is a disk-convex k\"ahler manifold) admits a meromorphic
extension to $\Omega_i$.
\end{theo}

As these results are already known for restrictions
of holomorphic functions defined in a neighborhood of $M$,
one natural idea is to apply the analytic disc techniques
in order to extend continuous CR functions on $M$ to a one sided neighborhood
of $M$. Then by deforming $\Omega_i$ (resp. $M$) 
in this one-sided neighborhood,
we are reduced to the case of holomorphic functions in the
neighborhood of $\overline \Omega_i$ (resp. $M$). 
This idea has already been applied
by J\"oricke, Merker or Porten 
in order to obtain many results about extension and removability
of singularities of CR functions. In the case of the study
of the Hartogs-Bochner phenomenon,  J\"oricke \cite{Joricke} proved
that compact hypersurfaces of $\mathbb{C}^n$ are globally minimal
(i.e. consist of a single CR orbit). Thus, using the propagation
results of Tr\'epreau \cite{trepreau} of analytic extension along
CR orbits, one obtains that CR functions defined on $M$
extends  holomorphically to a one sided neighborhood of $M$.
Thus, in the case of $\mathbb{C}^n$, the Hartogs-Bochner extension
theorem can be reduced to the classical Hartogs extension theorem
(this has been used for example in \cite{porten,sarkis} in order to prove
CR-meromorphic extension results). 
In the case of a compact Hypersurface of $P_n(\mathbb{C})$,
it is conjectured in \cite{sarkis} and also in \cite{dwil}
that compact hypersurfaces are also globally minimal but unfortunately
this is not known and is related to the following question of E. Ghys
(see \cite{ghys}):
{\it Does there exist a non-trivial
compact set laminated by Riemann surfaces in $P_2(\mathbb{C})$ ?}.\par\noindent
Indeed, in the case of a connected compact hypersurface $M$ of 
$P_n(\mathbb{C})$, CR orbits are either open subset of 
$M$ or injectively immersed
complex hypersurfaces whose closure is a compact subset
of $M$ laminated by complex manifolds of dimension $(n-1)$. 
Of course, if there exists no such laminated compact set in $P_n(\mathbb{C})$,
then $M$ has to be globally minimal (i.e. has only one open CR orbit).
Let $K$ be the union of all non open CR orbits of $M$. Then
$K$ is a laminated compact subset of $P_n(\mathbb{C})$. 
If $K=\emptyset$, $M$ is globally minimal and we are reduced to 
the result of \cite{sarkis}. If $K \neq \emptyset$, it is known that
$P_n(\mathbb{C}) \backslash K$ is Stein. Then, we apply the
boundary problem result of Chirka \cite{chirka} 
to the graph of CR functions over $M \backslash
K$ in order to obtain the needed holomorphic extension. 
In the case of CR maps having values
in a disc convex k\"ahler manifold, we follow the same idea
applying the boundary problem result
given in \cite{sarkis}.\par
\medskip
I would like to thank the referee and E. Chirka
for their remarks on this paper.
\section{Preliminaries}
\subsection{Decomposition into CR-orbits}

Let $M$ be an oriented and compact real hypersurface
of class $C^2$ of a complex manifold of dimension $n$.
For any point $p\in M$, we call $H_p(M)=T_p(M) \cap iT_p(M)$
the {\it holomorphic tangent space to $M$} at the point $p$ (where
$T_p(M)$ is the tangent space to $M$ at $p$).
As $M$ is of class $C^2$, the set
of holomorphic tangent spaces to $M$ is a vector bundle of
complex rank $n-1$. 
A $C^1$ curve $\gamma:[0,1] \rightarrow M$ is called
a {\it CR curve} if for any point $t \in [0,1]$, $\gamma'(t) \in H_{\gamma(t)}(M)$.
Let $x \in M$,  the set of points $y \in M$ which can be joined
to $x$ by a piecewise CR-curve is called 
the {\it CR-orbit ${\cal O}_{CR}(x)$} of $x$ in $M$.
It is well known that CR-orbits are CR-submanifolds injectively immersed
in $M$ and of the same CR dimension. Thus, for any point $x \in M$,
${\cal O}_{CR}(x)$ is either an open set (and we will say that
$M$ is globally minimal at the point $x$) or a  
complex manifold $\eta_x$ of dimension $(n-1)$ 
injectively immersed in $M$. 
In this last case, the CR orbits are tangent to the bundle $H(M)$
of complex tangent vectors to $M$.
As $M$ is of class $C^2$, $H(M)$ is of class $C^1$, thus any point $p \in M$
has a neighborhood $U_p$ such that $\eta_x$ is a product of the unit ball
of $\mathbb{C}^{n-1}$  by a topological set $T \subset \mathbb{C}$.
More precisely, $\overline \eta_x$ is a 
compact set laminated by complex manifolds of dimension $(n-1)$ (see \cite{ghys})~:
\begin{defi}
Let $N$ be a compact topological space
 and $\{U_i\}_{i \in I}$ be an open covering of $N$ 
such that for any $i \in I$, there exists an homeomorphism $h_i$ of $U_i$ on $\mathbb{B} \times T_i$
where $\mathbb{B}$ is the unit ball of $\mathbb{C}^n$ and $T_i$ is a topological space.
We say that these open sets define an atlas of a structure of lamination by complex
manifolds of dimension $n$ if the change of charts $h_{ij}=h_j \circ h_i^{-1}$ on their domain
of definition are of the following form:
$$h_{ij}(z,t)=(f_{ij}(z,t),\gamma_{ij}(t)),$$
where $f_{i,j}$ depends holomorphically of the variable $z$ and continuously of the variable $t$.
Two atlas on $N$ are equivalent if there union is an atlas.
A topological space is called laminated if an equivalence class
of atlases is given.
\end{defi}
\
Moreover, let $$K=\{ x \in M; M \mbox{ is not globally minimal at the point } x\}$$
then $K$ is also a compact set of $M$ (as its complement is open by definition) 
and is laminated by complex manifolds of dimension $(n-1)$ (see \cite{sussman,jori} for
a precise study of the structure of CR orbits).\par
As remarked in \cite{Joricke}, in order to prove the global minimality
of compact hypersurfaces of $\mathbb{C}^m$, one has to show that
there exists no such laminated compact sets in Stein manifolds~:

\begin{prop}
Let $X$ be a Stein manifold. Then there is no laminated
compact set $Y \subset X$.
\end{prop}
\begin{demo}
Let us suppose that there exists such a compact set in $X$.
By embedding $X$ in $\mathbb{C}^n$, we obtain
a laminated compact set $Y \subset \mathbb{C}^n$.
Let $r>0$ be the infimum of the reals $s>0$ such
that $Y \subset \overline {B(0,s)}$ where
$B(0,s)$ is the ball of center the origin and radius $s$.
Let $z \in Y \cap \partial B(0,r)$, let $\mathbb{C}_z$ be the
complex line containing the segment $[0,z]$ and let 
$\pi: \mathbb{C}^n  \rightarrow \mathbb{C}_z$ the 
projection on $\mathbb{C}_z$. Let $D_z$ be a complex manifold 
contained in $Y$ and passing through the point $z$. Then the restriction
of $\pi_z$ on $D_z$ is a non constant holomorphic function whose modulus
has a maximum at the point $z$, this contradicts the maximum principle.
\fdem
\end{demo}

\subsection{Laminated compact sets of $P_n(\mathbb{C})$}

Let $Y$ be a compact subset of $P_n(\mathbb{C})$
laminated by complex manifolds of dimension $(n-1)$.
Then by definition of $Y$, $P_n(\mathbb{C}) \backslash Y$
is pseudoconvex (as at any point of its boundary there
is a piece of complex hypersurface included in the boundary).
So, according to \cite{takeuchi, fujita,
  fujita2, kiselman}, $P_n(\mathbb{C}) \backslash Y$ is Stein.
As a direct consequence, we obtain~:
\begin{prop}
Let $Y$ be a compact set of $P_n(\mathbb{C})$ ($n \geq 2$)
laminated by complex manifolds of dimension $(n-1)$
and $f:Y \rightarrow \mathbb{C}$ be a continuous function
on $Y$ which restriction on any complex manifold contained
in $Y$ is holomorphic. Then $f$ is constant on $Y$.
\end{prop}
\begin{demo}
In order to prove that $f(Y)$ contains only one point,
it is sufficient to prove that its topological boundary $\partial f(Y)$
contains only one point. First, let us remark that
for any point $x \in \partial f(Y)$, $f^{-1}(x)$
is a laminated compact set of  $P_n(\mathbb{C})$.
Indeed, let $y \in Y$ be a point such that $f(y)=x$.
From the open mapping theorem, $f$ is constant
on the maximal complex manifold passing through the point $y$.
Thus, it is constant on its closure and we obtained
that $f^{-1}(x)$ is a laminated compact subset of $Y$  that we will
note $Y_x$.
Now, let suppose that $\partial f(Y)$ contains two different points
$x_1$ and $x_2$. Then the sets $Y_{x_1}$ and $Y_{x_2}$ 
are two compact sets of  $P_n(\mathbb{C})$ laminated by complex
manifolds of dimension $(n-1)$ which does not intersect
as $f$ takes different value on each one.
But,  $P_n(\mathbb{C}) \backslash Y_{x_1}$ is a pseudoconvex
open set of $P_n(\mathbb{C})$, following \cite{takeuchi, fujita,
  fujita2, kiselman},  $P_n(\mathbb{C}) \backslash Y_{x_1}$ is 
Stein and $Y_{x_2} \subset P_n(\mathbb{C}) \backslash Y_{x_1}$
which contradicts proposition $2$.
\fdem
\end{demo}

Let $M$ be a real hypersurface of a complex manifold $X$,
then for any point $x \in M$, there exists an open connected
neighborhood $V_x$ of $x$ such that $V_x \backslash M$ is
a disjoint union of two connected open sets (which are
called {\it one sided neighborhood of $M$ at the point $x$}).
We will say that $\cal W$ is a {\it one-sided neighborhood of $M$}
if for any point $x \in M$, $\cal W$ contains a one sided
neighborhood of $M$ at $x$ (the side can change).
Applying proposition $2$ and the results on propagation
of CR extension of \cite{trepreau, merker, Joricke},
we obtain the following proposition:

\begin{prop}
Let $M$ be a connected $C^2$ real hypersurface of 
$P_n(\mathbb{C})$ which divides $P_n(\mathbb{C})$
into two connected open sets. Let $K$ be the set of points
$x \in M$ such that $M$ is not globally minimal at the point $x$.
Then continuous CR functions on $M$ are constant on $K$.
Moreover, any continuous CR map defined on $M\backslash K$
and with values in any complex manifold extends holomorphically 
to a one sided neighborhood of any point of $M\backslash K$.
\end{prop}

\subsection{  Holomorphic decomposition of CR functions}

In this section, we give a proof (communicated to us
by C. Laurent-Thiebault) of the classical decomposition theorem
for CR functions as difference of boundary values of holomorphic functions:
\begin{prop}
Let $X$ be a complex manifold and $M \subset X$
be a $C^k$ compact real hypersurface which divides $X$ in
two open sets $\Omega_1$ and $\Omega_2$. Let $f:M \rightarrow \mathbb{C}$
be a CR function of class $C^{l+\alpha}$ ($0 \leq l < k$, $0<\alpha<1$) on $M$.
Let us suppose that $H^{0,1}(X)=0$. Then 
there exists two holomorphic functions $f_1$ and $f_2$ defined respectively
on $\Omega_1$ and $\Omega_2$ such that:
\begin{enumerate}
\item $\forall i \in \{1,2\}$, $f_i \in C^{l+\alpha}(\overline \Omega_i)$.
\item $f=f_1|_M-f_2|_M$.
\end{enumerate}
\end{prop}
\begin{demo}
Let $T=f[M]^{0,1}$ where $[M]^{0,1}$ is the part of
bidegree $(0,1)$ of the integration current over $M$.
As $f$ is a CR function, the current $T$ is 
$\overline \partial$-closed. The hypothesis $H^{0,1}(X)=0$
implies that there exists a distribution $S$ such that
$$\overline \partial S=T \mbox{ in } X.$$ The support of the current
$T$ being included in $M$, $\overline \partial S=0$
on $X \backslash M$. Thus $S$ defines two holomorphic functions
$f_1$ and $f_2$ defined respectively in $\Omega_1$ and $\Omega_2$.
If $z_0  \in M$, let us consider a neighborhood $V$ of $z_0$
biholomorphic to a ball in $\mathbb{C}^n$. One can then solve
the $\overline \partial$ on $V$ and obtain a distribution
$S_0$ on $V$ such that 
$$\overline \partial S_0=T \mbox{ in } V.$$
Thus, on $V$ we have that $$\overline \partial (S-S_0)=0$$
which implies that $S-S_0$ is a holomorphic function on $V$ 
and so is $C^\infty$. Thus, the regularity of $S$ in a neighborhood of $z_0$
is the same than the regularity of $S_0$ itself which jump 
is $C^{l+\alpha}$
as it can be checked using the kernels of Henkin in the ball. 
\fdem
\end{demo}

\subsection{Complex boundary problem}

Let $X$ be a complex Riemannian manifold of dimension $n$
and $M$ be a closed and oriented $C^1$ submanifold of $X$
of dimension $2p-1$ with $p \geq 1$
(we will note $[M]$ the integration current
associated to $M$). 
We will call {\it holomorphic $p$-chain}
any locally finite linear combination
of analytic subsets of $X \backslash M$ with integer coefficients. 
Of course, holomorphic
$p$-chains define closed currents of bidimension $(p,p)$
of $X \backslash M$. The volume of a holomorphic $p$-chain
$[T]=\sum n_i[T_i]$, is the expression
$$ \mbox{Vol }[T]=\sum |n_j|\mbox{Vol }V_j$$
where $\mbox{Vol }V_j$ (or ${\cal H}_{2p}(V_j)$) 
is the $2p$-dimensional Hausdorff measure of the analytic set $V_j$.
The volume  $\mbox{Vol }[T]$ is also equal to the mass of
the associated current $[T]$.
If a holomorphic $p$-chain is of locally finite mass
in $X$, the associated current in $X \backslash M$ will have an
extension to a current of $X$ which is not closed in general in $X$.
The question to find necessary and sufficient conditions
for $[M]$ to be the boundary (in the current sense) of a holomorphic
$p$-chain of $X \backslash M$ of locally finite volume in $X$
is called the  {\it complex boundary problem}.\par

Two necessary conditions for $[M]$ to have a solution
to the boundary problem are that $[M]$ 
can be decomposed into the sum of two currents 
of bidimension $(p,p-1)$ and $(p-1,p)$
(such a current will be said {\it maximally complex})
and that $[M]$ is a closed current.
Indeed, suppose that there exists a holomorphic 
$p$-chain $[T]$ of $X \backslash M$ such that 
$[T]$ is of locally finite volume in $X$ and verify $[M]=d[T]$.
Then necessarily we have
 $${\rm d}[M]={\rm d}({\rm d}[T])=0$$
and as $[T]$ is a current of bidimension $(p,p)$:
$$[M]={\rm d}[T]=(\partial+\overline \partial)[T]=\partial [T]+
\overline \partial [T]$$ is maximally complex.
Of course, the property for the current $[M]$ to be maximally complex 
is equivalent for the manifold $M$ to verify that
for any point $p \in M$, $\dim_{\mathbb{C}}H_p(M)=(p-1)$
where $H_p(M)$ is the holomorphic tangent bundle to $M$ at the point $p$.\par
In the case of $X=\mathbb{C}^n$, $p \geq 2$ and $M$ is compact, 
Harvey and Lawson \cite{harvey} proved that this two conditions
are in fact sufficient for the boundary problem for $M$ to have a
solution. Then, many authors studied the boundary problem
in more general manifolds (see for example \cite{chirka, henkin, dinh,
  sarkis2, sarkis3}). In this section we would like to mention
the following two results that will be used in the present article:

\begin{prop}[Chirka \cite{chirka}]
Let $Y$ be a polynomially convex compact set of $\mathbb{C}^n$
and $\Gamma$ a closed,  oriented and maximally complex $C^1$ submanifold
of $\mathbb{C}^n \backslash Y$ of dimension $2p-1$ ($p \geq 2$)
such that $\Gamma \cup Y$ is a compact set of $\mathbb{C}^n$. 
Then there exists a unique  
holomorphic $p$-chain $[T]$ of $\mathbb{C}^n \backslash (Y \cup M)$, 
of locally finite volume in $\mathbb{C}^n \backslash Y$ which is the solution
to the boundary problem for $[\Gamma]$ in $\mathbb{C}^n \backslash Y$
(i.e. $d[T]=[\Gamma]$).
\end{prop}

We will say that a complex manifold $Y$ is {\it disk-convex}
if for any compact set $L \subset Y$, there exists
a compact set $\widehat L$ such that,
for any irreducible analytic subset $S$ of $Y \backslash L$
verifying  $S \cup L$ is a compact subset of $Y$ and 
$\overline S \cap L \neq \emptyset$, we have $S \subset \widehat L$.
For example, any compact or holomorphically convex complex manifold
is disc-convex. \par

\begin{prop}[\cite{sarkis2}]
Let $U$ be a connected complex manifold of dimension $(n-1)$
and $\omega$ be a disc-convex k\"ahler manifold .
Let $\Gamma$ be a $C^\infty$ 
smooth maximally complex submanifold of $U \times \omega$
of real dimension $(2n-1)$. Let us note $\pi: U \times \omega
\rightarrow U$ the projection on the first member.
Let us suppose that $\Gamma$ verifies the following properties~:
\begin{enumerate}
\item [a.] For every compact $K \subset U$, $\Gamma \cap \pi^{-1}(K)$ is
a compact of $U \times \omega$.
\item [b.] For every $z \in U$, the set of points of
$\gamma_z= \Gamma \cap \{z\}\times \omega$   such that $\Gamma$ is
not transversal to $ \{z\}\times \omega$ is finite.
\item [c.] For every $z \in U$, $\gamma_z$ is a piecewise smooth curve
which have a finite set of singular points.
\end{enumerate}
Then, the following propositions are equivalent~:
\begin{enumerate}
\item There exists a non empty open set $O \subset U$ such that $[\Gamma]$
admits a solution to the boundary problem in $O \times \omega$
(i.e. there exists a holomorphic $n$-chain $[T]$ of
$(O \times \omega) \backslash \Gamma$, of locally finite mass in
$O \times \omega$, such that for any compact
$K \subset O$, $\overline T \cap \pi^{-1}(K)$ is a compact of
$O \times \omega$ and such that $d[T]= [\Gamma]$ (in the current sense).
\item There exists a closed subset  $F \subset U$ of $(2n-3)$
dimensional Hausdorff null measure such that for any point 
$z \in U \backslash F$, there exists an open set
$O_z \subset U$ such that $[\Gamma]$ admits a solution to the
boundary problem in $O_z \times \omega$.
\end{enumerate}
\end{prop}

The properties $b.$ and $c.$ of the proposition are generic in the
category of smooth manifolds of a product space. Thus any smooth manifold
has a deformation which verifies these properties. In particular,
applying this for graphs of holomorphic maps, we obtain the
following proposition~:

\begin{prop}
Let $M$ be a connected $C^2$ hypersurface of a complex
manifold $X$ of complex dimension $n$.
Let $f$ be a non constant 
holomorphic map having values in a complex manifold
$Y$ and defined
in a connected neighborhood of $M$.
Then there
exists a deformation $\widetilde M$ of $M$ in $X$
such that the graph 
$$\widetilde \Gamma_f= \{(w,z) \in Y \times X; z \in \widetilde M, w=f(z)\}$$ 
of the  restriction of $f$ on $\widetilde M$
verify the conditions $b.$ and $c.$ 
of the previous proposition.
\end{prop}

\section{ Proof of the main results}

Let us note
$$K=\{ x \in M; {\mathcal{O}}_{CR}(x) \mbox{ is not an open subset of } M\}.$$

\subsection{ Case $M$ is globally minimal (i.e. $K=\emptyset$)}

According to proposition $3$, if the compact set $K$ is empty,
continuous CR maps on $M$ extends holomorphically to
a one sided neighborhood of $M$.
Thus, by deforming $M$ into this one sided neighborhood and
by remarking that if $x \in M$ is such that CR maps
on $M$ extends to its two sides then they are restrictions
of holomorphic map in the neighborhood of $x$, we are reduced
to the following result:

\begin{prop}[\cite{sarkis,sarkis2}]
Let $M$ be a connected $C^1$ hypersurface of $P_n(\mathbb{C})$ ($n \geq 2$)
which divides $P_n(\mathbb{C})$ in two connected open sets $\Omega_1$ and $\Omega_2$.
Suppose that there exists a non constant holomorphic function defined
in a connected neighborhood of $M$.
Then holomorphic maps having values in disk-convex k\"ahler manifold
and defined in a connected
open neighborhood of $M$ extends meromorphically to $\Omega_1$ or 
to $\Omega_2$.
\end{prop}
\begin{demo}
We will give a 
proof of the proposition which is slightly
different of the one in \cite{sarkis,sarkis2}. 
Let $V$ be the connected open neighborhood of $M$
on which there exists a non constant holomorphic function
$g:V \rightarrow \mathbb{C}$.
\begin{lemme}
There exists a connected neighborhood $\widetilde V$ of $M$,
which is relatively compact in $V$ and two holomorphic
functions $f_1 \in \mathcal{O}(\Omega_1 \cup \widetilde V)$ 
and $f_2 \in \mathcal{O}(\Omega_2 \cup \widetilde V)$
such that $g|_{\widetilde V}=f_1|_{\widetilde V} - f_2|_{\widetilde V}$.
\end{lemme}
\begin{demo}
Let $\widetilde V$ be a connected open neighborhood
of $M$ which is relatively compact in $V$.
Let $\phi$ be a smooth function defined on $P_n(\mathbb{C})$
such that $\mbox{supp }\phi \subset V$ and $\phi|_{\widetilde V} \equiv 1$.
For any $i \in \{1,2\}$, 
let us consider the smooth forms $\omega_i$ defined by
$$\omega_i=\{
\begin{array}{c}
\bar \partial (g\phi) \mbox{ on } \Omega_i \\
0 \mbox{ on } P_n(\mathbb{C}) \backslash \Omega_i
\end{array}.
$$
Then, for each $i \in \{1,2\}$,
 $\omega_i$ is a $\bar \partial$-closed $(0,1)$ smooth form which support
is included in $\Omega_i \backslash \widetilde V$.
As $H^{0,1}(P_n(\mathbb{C}))=0$, there exists a smooth function
$u_i$ defined on $P_n(\mathbb{C})$ such that $\omega_i=\bar \partial u_i$.
Of course, we have $$\bar \partial (g\phi)=\omega_1+\omega_2
=\bar \partial u_1+ \bar \partial u_2=\bar \partial (u_1+u_2)$$
So the smooth function $u_1+u_2-g\phi$ is holomorphic on $P_n(\mathbb{C})$
and thus is constant. Let us note $c$ the constant
that verify $$c=u_1+u_2-g\phi.$$
Then the functions $f_1$ defined on $\Omega_1 \cup \widetilde V$
by $f_1=g\phi-u_1-c/2$ and $f_2$ defined on $\Omega_2 \cup \widetilde V$
by $f_2=u_2-g\phi-c/2$ 
are holomorphic on there respective domains of definition
and verify that on $\widetilde V$ (as $\phi|_{\widetilde V}\equiv 1$)
we have  $f_1-f_2= g$. \fdem
\end{demo}

As $g$ is supposed non constant, one of the two
holomorphic functions $f_1$ and $f_2$ has to be also non constant.
Let us suppose for example that $f_1$ is non constant.
Then, according to  \cite{takeuchi, fujita,
  fujita2, kiselman}, the envelope 
of holomorphy $W$ of $\Omega_1 \cup \widetilde V$ is Stein.
So the domain $\Omega_1 \cup \widetilde V$ embeds in its envelope
and in particular $\Omega_1$ can be seen as a bounded
domain of a Stein space $W$ and we are reduced to the classical
Hartogs-Levi extension theorem 
(see \cite{ivashkovitch} or \cite{sarkis} for 
a proof of the Hartogs-Levi extension theorem
in the case of meromorphic maps having values in disk-convex
k\"ahler manifolds).
\end{demo}

\subsection{ Case $M$ is not globally minimal (i.e $K \neq \emptyset$)}

First, let us  remark that the theorems $1$, $2$ and $3$ 
are trivial in the
case that there is no non-constant $C^1$ CR functions on $M$.
So in all the following we will always assume that
$g: M \rightarrow \mathbb{C}$ is a non-constant $C^1$ CR function.
\begin{lemme}
The compact set $K$ must verify the following properties:
\begin{enumerate}
\item $K$ is a compact set laminated by complex
manifolds of dimension $(n-1)$.
\item The CR function $g$ is constant on $K$ (we can suppose that $g(K)=0$).
\item  $K$ is of null $(2n-1)$-dimensional Hausdorff measure.
\item The open set $U=P_n(\mathbb{C}) \backslash K$ is Stein.
\end{enumerate}
\end{lemme}
\begin{demo} 
The two first points are consequence of proposition $3$.
According to proposition $4$, $g=f_1-f_2$ where
$f_i \in {\mathcal O}(\Omega_i) \cap C^0(\overline \Omega_i)$.
As $g$ is not constant then one of the two functions $f_i$
has also to be non constant. But, as they are constant
on $K$, the set $K$ has to be of null measure in $M$ (which proves
the third point). 
Finally, As the compact set $K$ is supposed non empty,
$$U=P_n(\mathbb{C})\backslash K$$ is a pseudoconvex open subset of
$P_n(\mathbb{C})$. According to \cite{takeuchi, fujita,
  fujita2, kiselman}, $U$ is Stein. 
\end{demo}

\subsubsection{Semi-local solution to the boundary problem}

Let us suppose that $f:(M\backslash K) \rightarrow \mathbb{C}$
is a continuous CR function. According to proposition $3$,
up to deforming $M$, we can always assume that $g$ and $f$
are smooth on $M \backslash K$.
We will prove that the graph of the restrictions of $f$
over the level sets $\{g=c\}$ admits solutions to the boundary problem.
More precisely, Let us consider the graph of the map $(g,f)$
over the set $(M \backslash \{g = 0\}) \subset M \backslash K$:
$$\widetilde \Gamma_{g,f}=
\{(w,y,z) \in (\mathbb{C} \backslash \{0\}) \times \mathbb{C} 
\times P_n(\mathbb{C}); z \in M\backslash \{g=0\}, w=g(z), y=f(z)\}$$ 
\begin{prop}
There exists a holomorphic $n$-chain $[\widetilde T]$
of  $((\mathbb{C}\backslash \{0\}) \times \mathbb{C} \times P_n(\mathbb{C}))
\backslash \widetilde \Gamma_{g,f}$, of locally finite mass  
in $(\mathbb{C}\backslash \{0\}) \times \mathbb{C} \times P_n(\mathbb{C})$ 
solution to the boundary problem for  $[\widetilde \Gamma_{g,f}]$
(moreover, for any compact $R \subset \mathbb{C}\backslash \{0\}$,
the set $\widetilde T \cap (R \times \mathbb{C} \times P_n(\mathbb{C}))$ 
is compact in $(\mathbb{C}\backslash
\{0\}) \times \mathbb{C} \times U$).\fdem
\end{prop}
\begin{demo}
We recall that $U=P_n(\mathbb{C}) \backslash K$ is Stein.
Let us note $D(0,\epsilon)$ the disc of center $0$ and radius $\epsilon$,
$C(0,\epsilon)$ its boundary and $\pi: \mathbb{C}
\times \mathbb{C} \times U \rightarrow \mathbb{C}$ be the projection on the first member.
For any $\epsilon >0$,  $\widetilde \Gamma_{g,f} 
\cap (C(0,\epsilon) \times \mathbb{C} \times U)$ is bounded, thus, there exists
a holomorphically convex compact set 
$B_\epsilon$ in $\mathbb{C} \times U$
such that $[\widetilde \Gamma_{g,f}
\cap (\mathbb{C} \backslash D(0,\epsilon))
\times \mathbb{C} \times U)] \cup D(0,\epsilon)\times B_\epsilon$
is a compact subset of $\mathbb{C} \times \mathbb{C} \times U$.
According to proposition $5$ with $Y=D(0,\epsilon) \times B_\epsilon$,
the boundary problem for $[\widetilde \Gamma_{g,f}]$
admits a unique solution in 
$(\mathbb{C}
 \times \mathbb{C} \times U) \backslash (Y \cup \widetilde
 \Gamma_{f})$.
By uniqueness of the solution and by leting $\epsilon$ tend
to zero, we obtain that the boundary problem for 
$[\widetilde \Gamma_{g,f}]$ has a unique solution 
$[\widetilde T]$ to the boundary problem
in $(\mathbb{C}\backslash \{0\}) \times \mathbb{C} \times U$.
As $[\widetilde \Gamma_{g,f}]$ has a solution to the boundary
problem in $(\mathbb{C} \backslash \{0\}) \times \mathbb{C}
\times U$ and as $U \subset P_n(\mathbb{C})$ we obtained
in fact that $[\widetilde T]$ is
a solution to the boundary problem for 
$[\widetilde \Gamma_{g,f}]$ in 
$(\mathbb{C} \backslash \{0\}) \times \mathbb{C} \times P_n(\mathbb{C})$.
\fdem
\end{demo}

For any $c \in \mathbb{C}\backslash \{0\}$,
let $\gamma_c=(\{c\}\times \mathbb{C} \times P_n(\mathbb{C}))
\cap \Gamma_{g,f}$ be the graph of $f$ over the level
set $\{g=c\}$. According to Sard's theorem, 
for almost all $c \in \mathbb{C} \backslash \{0\}$, 
$\gamma_c$ is a smooth curve and the 
intersection current $[\gamma_c]$ obtained by slicing $[\Gamma_{g,f}]$ with
the fiber $\{c\}\times \mathbb{C} \times P_n(\mathbb{C})$
is well defined and corresponds to the integration current
over $\gamma_c$. Moreover, for almost all
$c \in \mathbb{C} \backslash \{0\}$, the boundary in the current
sense of the intersection current
 (noted $[S_c]$) obtained by slicing the current
 $[\widetilde T]$ by the fiber
$\{c\}\times \mathbb{C} \times P_n(\mathbb{C})$
is equal to the intersection current of 
$[\widetilde \Gamma_{g,f}]$
by this same fiber. So we obtain:
\begin{coro}
For almost all $c \in \mathbb{C} \backslash \{0\}$, 
there exists a holomorphic $1$-chain $[S_c]$
of $(\{c\}\times \mathbb{C} \times U)
\backslash \gamma_{c}$, of finite mass in 
$\{c\} \times \mathbb{C} \times U$ solution to the boundary problem
for $[\gamma_c]$ (i.e. $d[S_c]=[\gamma_c]$).
\end{coro}

\subsubsection{Global solution to the boundary problem}

Let $g$ and $f$ be two
$C^1$ CR functions defined on $M$ and let us suppose
that $g$ is not constant.
Let $$\Gamma_{g,f}=\{(w,y,z) \in \mathbb{C} \times \mathbb{C}
\times P_n(\mathbb{C}); z \in M, w=g(z), y=f(z)\}$$ be the graph of
the map $(g,f)$. 
Let us prove that $[\Gamma_{g,f}]$, the integration current over $\Gamma_{g,f}$,
admits a solution to the boundary problem in 
$\mathbb{C} \times \mathbb{C} \times P_n(\mathbb{C})$. \par

\begin{prop} Let $g,f$ and $\Gamma_{g,f}$ as defined above.
Then there exists
an irreducible holomorphic $n$-chain $[T]$ 
of $(\mathbb{C} \times \mathbb{C}\times  P_n(\mathbb{C}))\backslash \Gamma_{g,f}$,
of finite mass in $\mathbb{C} \times \mathbb{C} \times \times P_n(\mathbb{C})$
which is solution to the boundary problem for $[\Gamma_{g,f}]$
(i.e. $[\Gamma_{g,f}]=d[T]$ in the current sense).
\end{prop}
\begin{demo}
According to proposition $9$, 
there exists a holomorphic $n$-chain $[\widetilde T]$, 
of locally finite mass  
of $(\mathbb{C}\backslash \{0\}) \times \mathbb{C} \times P_n(\mathbb{C})$ 
solution to the boundary problem for  $[\widetilde \Gamma_{g,f}]$.
Let us prove that $[\widetilde T]$ extends
to a solution to the boundary problem for 
$[\Gamma_{g,f}]$ in all $\mathbb{C} \times \mathbb{C} \times P_n(\mathbb{C})$.
For this, according to the Bishop's extension theorem (see \cite{bishop}),
it is enough to prove that $[\widetilde T]$ 
is of finite mass in $\mathbb{C} \times \mathbb{C} \times P_n(\mathbb{C})$
in order for it to extend to a holomorphic $n$-chain
of $(\mathbb{C} \times \mathbb{C} \times P_n(\mathbb{C})) 
\backslash \Gamma_{g,f}$.

\begin{lemme} The holomorphic $n$-chain
$[\widetilde T]$ is of finite mass in $\mathbb{C} 
\times \mathbb{C} \times P_n(\mathbb{C})$.
\end{lemme}
\begin{demo}
Let $\Omega$ be the k\"ahler form associated
to $\mathbb{C} \times P_n(\mathbb{C})$. Let $\{\phi_\epsilon\}$  be a family
 of smooth functions defined on $\mathbb{C}$ and verifying:
\begin{enumerate}
\item $0 \leq \phi_\epsilon \leq 1$
\item $\forall z,  |z| > 2\epsilon$ we have $\phi_\epsilon(z)=1$
\item $\forall z,  |z| < \epsilon$ we have  $\phi_\epsilon(z)=0$
\end{enumerate}
As $[\widetilde T]$ is a holomorphic $n$-chain, we have
$[\widetilde T]=\sum_{j \in J} n_j [T_j]$ where $[T_j]$
are integration current over analytic subsets of dimension $n$
of $((\mathbb{C} \backslash \{0\}) \times \mathbb{C} \times P_n(\mathbb{C})) 
\backslash \Gamma_{g,f}$.
The mass  $Vol([\widetilde T])$ of the current $[\widetilde T]$
is equal to the mass of the positive current
$\sum_{j\in J} |n_j|[T_J]$. 
So we have
$$Vol([\widetilde T])=\sum_{j \in J} |n_j| Vol(T_j)=
\sum_{j \in J} \frac{|n_j|}{n!} \int_{T_j} (\Omega+dz\wedge d\bar
z)^n$$
$$=\lim_{\epsilon \rightarrow 0} \sum_{j \in J} \frac{|n_j|}{n!}
\int_{T_j} \phi_\epsilon(\Omega +dz \wedge d\bar z)^n$$ 
For all $\epsilon > 0$, lets note
$\gamma_\epsilon=(\int_d^z \phi_\epsilon(u) du)\wedge d\bar z$ 
a primitive of the closed form $\phi_\epsilon dz \wedge d\bar z$
which vanishes in a neighborhood of $d$. 
We have
 $$\int_{T_j} \phi_\epsilon (\Omega+dz\wedge d\bar z)^n
= \int_{T_j} \phi_\epsilon \Omega^{n-1}\wedge dz\wedge d\bar z$$
(because $dz \wedge d\bar z$ is of maximal degree on $\mathbb{C}$
and $T_j$ has no vertical components of dimension $n$
so $\int_{T_j} \Omega^{n}=0$) 
$$=\int_{T_j} \Omega^{n-1} \wedge d\gamma_\epsilon 
= \int_{T_j}d( \Omega^{n-1}\wedge \gamma_\epsilon) $$
(because $\Omega$ is closed)
$$= \int_{\partial T_j} \Omega^{n-1} \wedge \gamma_\epsilon$$
Let us note $\{\Gamma_l\}_{l\in L}$ 
the connected components of 
$\widetilde \Gamma_{g,f}$.
Then we have (by grouping the connected components of the boundaries
of the $T_j$):
$$\sum_{j \in J} |n_j|\int_{\partial T_j} \Omega^{n-1} 
\wedge \gamma_\epsilon
\leq \sum_{l \in L} \pm \int_{\Gamma_l} \Omega \wedge \gamma_\epsilon
\leq \int_{\Gamma_{g,f}}|\Omega^{n-1}\wedge \gamma_\epsilon|$$
Since $\Gamma_{g,f}$ is of finite volume and, since the forms
$\Omega^{n-1}\wedge \gamma_\epsilon$ are uniformly bounded, we deduce
that $Vol([\widetilde T]) < \infty$ which
ends the proof of the lemma.\fdem
\end{demo}

Let us note $[T]$ the holomorphic $n$-chain
of $(\mathbb{C} \times \mathbb{C} \times P_n(\mathbb{C})) \backslash
\Gamma_f$ which is the simple extension of $[\widetilde T]$.
Let us call $$[B]=d[T]-[\Gamma_{g,f}].$$ 
Then $[B]$ is a flat current of dimension $(2n-1)$ 
which support is included in $K$ (which is of 
null $(2n-1)$-Hausdorff measure). So $[B]=0$ (see \cite{federer}
theorem 4.1.20 p.378) and thus
$[T]$ is a solution to the boundary problem for 
$[\Gamma_{g,f}]$ (i.e. $[\Gamma_{g,f}]=d[T]$ in 
$\mathbb{C} \times \mathbb{C} \times P_n(\mathbb{C})$). 
As $\Gamma_{g,f}$ is connected,
$[T]=\pm[A]$ where $A$ is the integration current over some
irreducible analytic subset $A$ of finite volume
of $(\mathbb{C} \times \mathbb{C} \times
P_n(\mathbb{C})) \backslash \Gamma_{g,f}$
(we will say that $[T]$ is an irreducible holomorphic chain
with multiplicity $\pm 1$).
\fdem
\end{demo}

\subsubsection{Hartogs-Bochner phenomenon in the smooth case}

In this section we will prove the following smooth
version of theorem $1$:
\begin{prop}
Let $\Omega \subset P_n(\mathbb{C})$ ($n\geq 2$) be a domain
with $C^2$ and connected boundary $M$. Suppose
there exists a non-constant holomorphic function 
$g \in \mathcal{O}(\Omega) \cap C^1(\overline \Omega)$
then $C^1$ CR function $f:M \rightarrow \mathbb{C}$ 
extends holomorphically to $\Omega$.
\end{prop}
\begin{demo}
According to proposition, there exists
a solution $[T]$ to the boundary problem for $[\Gamma_{g,f}]$
in $\mathbb{C} \times \mathbb{C} \times P_n(\mathbb{C})
\backslash \Gamma_{g,f}$.
Let us prove that $[T]$ is the graph of a holomorphic
extension of the map $(g,f)$ on $\Omega$.
Let us define
$$\Gamma_{g}=\{(w,z) \in \mathbb{C} \times P_n(\mathbb{C});
z \in M, w=g(z)\}$$
$$T_g=\{(w,z) \in \mathbb{C} \times P_n(\mathbb{C}); z \in
\Omega, w=g(z)\}.$$
Let $[\Gamma_{g}]$ and $[T_g]$ be the associated integration currents.
Of course, we have $d[T_g]=[\Gamma_{g}]$ (the orientation 
of $M$ is supposed chosen such that $M$ is the oriented boundary of 
$\Omega$).
Let $\Pi: \mathbb{C} \times \mathbb{C} \times P_n(\mathbb{C})
\rightarrow \mathbb{C} \times P_n(\mathbb{C})$ be the projection
defined by $\Pi(w,y,z)=(w,z)$. Then we have
$$d\Pi_*([T])=\Pi_*(d[T])=\Pi_*([\Gamma_{g,f}])=[\Gamma_{g}]=d[T_g].$$
\begin{lemme}
We have $\Pi_*([T])=[T_g]$. 
\end{lemme}
\begin{demo}
Let $[A]=\Pi_*([T])-[T_g]$, then $[A]$ is a
holomorphic $n$-chain of $(\mathbb{C} \times P_n(\mathbb{C}))
\backslash \Gamma_{g}$ which is closed in the current
sense in $\mathbb{C} \times P_n(\mathbb{C})$.
So $[A]$ defines a closed holomorphic $n$-chain
of $\mathbb{C} \times P_n(\mathbb{C})$ which support
$A$ is a compact analytic subset of $\mathbb{C} \times P_n(\mathbb{C})$.
The projection of $A$ on $\mathbb{C}$ is a compact analytic subset
of $\mathbb{C}$ and so is a finite number of points. So,
as $g$ is not constant, its projection is empty
and so is $A$. So $[A]=0$ and we proved the lemma. \fdem
\end{demo}

So as $[T]$ and $[T_g]$ are irreducible holomorphic $n$-chains
of multiplicity $+1$, the restriction $\Pi|_T$ defines a (weak) biholomorphism
between $T_g$ and $T$.
Let $\Pi_2: \mathbb{C} \times P_n(\mathbb{C}) \rightarrow P_n(\mathbb{C})$
be the projection on the second member. Then by definition of $T_g$,
the restriction $\Pi_2|_{T_g}$ is a biholomorphisme between
$T_g$ and $\Omega$. Finally, the restriction $(\Pi_2 \circ \Pi_1)|_{T}$
defines a biholomorphisme between $T$ and $\Omega$. Thus it
defines the graph of the wanted extension of the map $(g,f)$
and we have proved that $f$ extends holomorphically to $\Omega$.\fdem
\end{demo}

Of course, if there exists a $C^{1+\alpha}$ ($0<\alpha<1$)
CR function $g$ defined
on $M$, according to proposition $4$, $g=f_1-f_2$
where $f_i$ are holomorphic functions defined
on each side of $M$ and which are $C^1$ up to the boundary.
As $g$ is supposed non constant, one of these two functions
has also to be non constant and we obtain the following corollary:
\begin{coro}
Let $M$ be a connected $C^2$ hypersurface of $P_n(\mathbb{C})$ ($n\geq 2$)
which divides $P_n(\mathbb{C})$ in two connected open sets $\Omega_1$
and $\Omega_2$. Suppose that there exists a $C^{1+\alpha}$ non-constant
CR function $g$ defined on $M$. 
Then there exists $i \in \{1,2\}$ 
such that $C^1$ CR functions on $M$ extends holomorphically to $\Omega_i$.
\end{coro}

\subsubsection{ Connectness of $M \backslash K$}

In the proof of the smooth case, one important point 
that shows that the solution of the boundary problem gives
an univalued holomorphic extension is that,
as $M$ is connected, $\Gamma_{g,f}$ is also connected
and so its solution to the boundary problem is irreducible.
In fact, the proof of the smooth case
implies implicitly that $M\backslash K$ is connected:

\begin{prop}
Let $M$ be a $C^2$ connected hypersurface of $P_n(\mathbb{C})$ ($n\geq 2$)
which divides $P_n(\mathbb{C})$ in two connected open sets
$\Omega_1$ and $\Omega_2$. Suppose that there exists
a non-constant $C^1$ CR function $g$ defined on $M$.
Let us suppose that $g(K)=\{0\}$. 
Then the subset $M \backslash \{g=0\}$ is connected.
In particular the set $M \backslash K$ is connected.
\end{prop}
\begin{demo}
According to proposition $3$, $g$ extends
holomorphically in a one sided neighborhood of
any point of $M \backslash K$. So the level set
$\{g=0\}$ cannot disconnect any connected
component of $M \backslash K$.
Let $M_2$ be a connected component of $M \backslash \{g=0\}$
on which $g$ is not constant and let suppose by contradiction
that $M_2$ is not equal to all of $M \backslash \{g=0\}$.
So, the open set  $V=P_n(\mathbb{C})\backslash \overline M_2$ is connected.
We can apply the proof of proposition $10$ to the graph
$\Gamma_2$ of $g$ over $M_2$ (take for example $f \equiv 0$). 
So we find an irreducible
holomorphic $n$-chain $[T_2]$ of $(\mathbb{C}
\times P_n(\mathbb{C})) \backslash \overline \Gamma_2$ of finite mass in
$\mathbb{C} \times P_n(\mathbb{C})$
and such that $d[T_2]=[\Gamma_2]$ in $(\mathbb{C} \backslash
\{0\})\times P_n(\mathbb{C})$. 
Let $\Pi: \mathbb{C} \times P_n(\mathbb{C}) \rightarrow P_n(\mathbb{C})$
be the projection on the second member. Let $$[C]=\Pi_*([T_2]),$$
then $[C]$ is a holomorphic $n$-chain of
$V$ which is solution to the boundary problem for $[M_2]$
in $P_n(\mathbb{C}) \backslash \{g=0\}$  as we have:
$$d[C]=d\Pi_*([T_2])=\Pi_*(d[T_2])=\Pi_*([\Gamma_2])=[M_2].$$
But as the open set $V$ is supposed
connected, we have that $[C]=k[V]$ where $k \in \mathbb{Z}$ and $[V]$
is the integration current over $V$. But, in this case
the current $d[C]=kd[V]$ has its supports in the set $\{g=0\}$ which contradicts
the fact that $d[C]=[M_2]$ in $P_n(\mathbb{C})\backslash \{g=0\}$.\fdem
\end{demo}

\subsubsection{Proof of theorem $2$: continuous case}

According to proposition $4$, $g=f_1|_M-f_2|_M$
where $f_i \in \mathcal{O}(\Omega_i) \cap C^0(\Omega_i)$.
As $g$ is not constant, one of these two functions
is also non constant (let suppose for example that
$f_1$ is not constant, that $f_1(K)=\{0\}$  and that the orientation
of $M$ has been chosen such that $M$ is the oriented boundary
of $\Omega_1$).

\begin{prop}
Any continuous CR function $f:(M \backslash K) \rightarrow \mathbb{C}$
extends holomorphically to $\Omega_1 \backslash \{f_1=0\}$.
Moreover, if $f$ is bounded, then its extension is also bounded and
according to the Riemann extension theorem, it extends
holomorphically to all of $\Omega_1$.
\end{prop}
\begin{demo}
According to proposition $3$, we can assume that
all the considered functions $f$ and $f_1$ are smooth in $M \backslash K$.
Let us define 
$$\widetilde \Gamma_{f_1,f}=\{(w,y,z) \in \mathbb{C} \times 
\mathbb{C} \times P_n(\mathbb{C}); z \in M \backslash \{f_1=0\},
w=f_1(z), y=f(z)\}.$$
Then according to proposition $9$, there
exists a holomorphic $n$-chain $[\widetilde T]$
of locally finite mass in $(\mathbb{C} \backslash \{0\})
\times \mathbb{C} \times P_n(\mathbb{C})$ such
that $d[\widetilde T]=[\widetilde \Gamma_{f_1,f}]$.
According to proposition $12$, as $g$ is of class $C^1$ and non-constant,
$M \backslash K$ is connected. So $M \backslash \{f_1=0\}$ is also
connected and so is $\widetilde \Gamma_{f_1,f}$.
So the holomorphic $n$-chain $[\widetilde T]$ has to be irreducible
and with multiplicity $\pm 1$.
For almost all $c \in \mathbb{C} \backslash \{0\}$, 
the intersection current $[\gamma_c]$ of $[\widetilde \Gamma_{f_1,f}]$
by the fiber $\{c\} \times \mathbb{C}\ \times P_n(\mathbb{C})$ 
is well defined.  According to corollary $1$, for almost all 
$c \in \mathbb{C}\backslash \{0\}$, we have that $[\gamma_c]=d[S_c]$
 where $[S_c]$ is the intersection current of $[\widetilde T]$
by the fiber $\{c\}\times \mathbb{C} \times P_n(\mathbb{C})$.
Let $\Pi_3: \mathbb{C} \times \mathbb{C} \times P_n(\mathbb{C})
\rightarrow P_n(\mathbb{C})$ be the projection on the last member.
For almost all $c \in \mathbb{C} \backslash \{0\}$ 
the integration currents $[\gamma_c^1]$ and $[S_c^1]$
respectively over the level sets $\{f_1|_M=c\}$ and  $\{f_1=c\}$
are well defined and verify that $[\gamma_c^1]=d[S_c^1]$.
Then, we have for almost all $c \in \mathbb{C} \backslash \{0\}$,
$$d{\Pi_3}_*([S_c])={\Pi_3}_*(d[S_c])={\Pi_3}_*([\gamma_c])
=[\gamma_c^1]=d[S_c^1]$$

By the uniqueness of the solution to the boundary problem
in the Stein manifold $\{c\}\times \mathbb{C} \times U$
we have that
$${\Pi_3}_*([S_c])=[S_c^1].$$
As the current $[\widetilde T]$ is of multiplicity $\pm 1$,
the currents $[S_c]$ have also to be of multiplicity $\pm 1$.
But, as by construction, the current $[S_c^1]$ is of multiplicity
$+1$ and verify that ${\Pi_3}_*([S_c])=[S_c^1]$, the current $[S_c]$
has also to be of multiplicity $+1$.
This proves that $f$ extends holomorphically
on almost all the level lines $S_c^1$ of $f_1$ which proves
(by continuity of the solution to the boundary problem)
that $f$ extends holomorphically to $\Omega_1 \backslash \{f_1=0\}$.\\
In the case that  $f$ is bounded, 
by construction of the solution to the boundary problem
the extension has also to be bounded.
\fdem
\end{demo}

The only remaining point in the proof of theorem $2$
is the regularity of the extension.
According to \cite{harvey} theorem 5.2 p.249,
the regularity up to the boundary of the holomorphic extension
is the same than the one of the considered CR function on $M$.
This ends the proof of the assertions $1$ and $2$
of theorem $1$. Concerning the third assertion,
let suppose that $f$ is a holomorphic function
on $\Omega_j$ which is continuous on $\overline \Omega_j$.
Then, $f|_M$ admits a holomorphic extension
to $\Omega_i$ which is continuous up to the boundary. 
So  $f$ extends holomorphically to all of $P_n(\mathbb{C})$
and consequently have to be constant.

\subsubsection{Case of CR maps: proof of theorem $3$}

Let us suppose moreover that $n=2$.
Let $f$ be
a $C^1$ CR map $f:M \rightarrow X$ where
$X$ is a disk-convex K\"ahler manifold. Let us consider the graph
of the map $(g,f)$:
$$\Gamma_{g,f}=\{(w,y,z) \in \mathbb{C} \times X \times P_2(\mathbb{C}); 
z \in M, w=g(z), y=f(z)\}$$
and let note $[\Gamma_{g,f}]$ the integration current
associated to it. Let us note  $\widetilde \Gamma_{g,f}$
the restriction of  $\Gamma_{g,f}$ to  
$(\mathbb{C}\backslash \{0\}) \times X \times U$
and  $[\widetilde \Gamma_{g,f}]$ the associated current.\par
According to theorem $2$, their exists $i \in \{1,2\}$
such that the CR function $g$ admits 
a holomorphic extension to $\Omega_i$ which is $C^1$ up to the boundary.
So, all the proof of section $3.2.3$ can be applied
except the use of Chirka's solution to the boundary 
problem (proposition $5$) 
which has to be replaced by the use of proposition $6$.
So, we only have to prove the following proposition:

\begin{prop}
There exists holomorphic $p$-chain $[\widetilde T]$ (which support is
noted $\widetilde T$), of locally finite mass  of
$(\mathbb{C}\backslash \{0\}) \times X \times U$ solution
to the boundary problem for  $[\widetilde \Gamma_{g,f}]$
(i.e. for any compact $Y \subset \mathbb{C}\backslash \{0\}$,
the set $\widetilde T \cap Y \times U$ is compact in $\mathbb{C}\backslash
\{0\} \times X \times U$ and $d[\widetilde T]=[\widetilde \Gamma_{g,f}]$
in the current sense).
\end{prop}
\begin{demo}
Let $\pi: \mathbb{C} \times X
\times U \rightarrow \mathbb{C}$ be the projection on the first member.
We will solve the boundary problem for $[\widetilde \Gamma_{g,f}]$
by applying the proposition $6$ to $[\widetilde \Gamma_{g,f}]$
with $U=\mathbb{C}$ and $\omega=\mathbb{C} \times X$.
As $g$ is constant on $K$, the hypothesis $a.$ of
proposition $6$ is verified. According to proposition $3$,
the map $(g,f)$ extends holomorphically to a one sided neighborhood
of $M \backslash K$. So, up to deforming $M$ in this one sided neighborhood
one can always suppose that the hypothesis $b.$ and $c.$
are verified. The function $g$ being continuous on $M$, 
it is bounded (let $N$ the maximum of its modulus) and
let $O$ be the complement in $\mathbb{C}$ of the closed disc
$\overline{D(0, 2N)}$
of center $0$ and radius $2N$. 
Then $\widetilde \Gamma_{g,f} \cap \pi^{-1}(O)=\emptyset$
and thus $[\widetilde \Gamma_{g,f}]$ have the null current as solution to
the boundary problem in $\pi^{-1}(O)$. According to 
proposition $6$,
there exists a closed set $F \subset \mathbb{C} \backslash \{0\}$
of Hausdorff $1$-dimensional null measure such that any point
$z \in \mathbb{C}\backslash (\{0\} \cup F)$, have a neighborhood $V_z$
such that $\widetilde \Gamma_{g,f}$ has a solution to the boundary problem
in $\pi^{-1}(V_z)$. \par
\
\begin{lemme}
$[\widetilde \Gamma_{g,f}]$ has a unique solution $[\widetilde T]$ 
to the boundary problem in \\ $\pi^{-1}(\mathbb{C}\backslash (\{0\} \cup F))$.
\end{lemme}
\begin{demo}
According to \cite{sarkis2} lemma 2.8, there exists a maximal open connected
set $V_{max} \subset \mathbb{C} \backslash \{0\}$ with
$V_{max} \supset O$ and  such that $[\widetilde \Gamma_{g,f}]$
has a solution $[\widetilde T]$ to the boundary problem in $\pi^{-1}(V_{max})$.
Let us prove by contradiction that 
$V_{max}= \mathbb{C}\backslash (\{0\} \cup F)$. Indeed,
suppose there exists a point $z$ in $\partial V_{max}
\cap \mathbb{C}\backslash (\{0\} \cup F)$. Then by
definition of $F$, there exists an open and connected neighborhood 
$V_z$ of $z$ such that $[\widetilde \Gamma_{g,f}]$ admits a solution $[R]$
to the boundary problem in $\pi^{-1}(V_z)$. Because of the
maximality assumption on $V_{max}$ there exists a connected
component $V$ of $V_{max} \cap V_z$ such that $[R]$ and 
$[\widetilde T]$ do not coincides on $\pi^{-1}(V)$. The current
$[L]=[R]-[\widetilde T]$ is a closed rectifiable $(2,2)$-current.
Thus, according to Bishop-King-Harvey-Shiffman-Alexander
\cite{bishop,king,harvey2,shiffman,alexander},
it is a holomorphic $2$-chain of
$\pi^{-1}(V)$. Thus $[L]=\sum_{j \in J} [L_j]$
where $L_i$ are analytic subsets of dimension $2$ of $\pi^{-1}(V)$.
By construction, for any compact $Y \subset V$, $\pi^{-1}(Y) \cap L_j$
is compact in $\pi^{-1}(V)$. Let $j \in J$, $\Gamma$ be a connected component
of $\Gamma_{g,f} \cap L_j$ and $w \in V$ be a point such that
$\gamma=\pi^{-1}(\{w\}) \cap \Gamma$ is a non empty smooth curve.
Then $S=\pi^{-1}(w) \cap L_j$ is a compact curve of $\{w\} \times X \times U$
containing $\gamma$.
According to the proper map theorem, the projection of $S$ on $U$ 
is a compact analytic subset of $U$. As $U$ is a stein manifold, this
projection has to be a point but it should contain the projection
of $\gamma$ which is not a point as the projection on $U$ of 
$\widetilde \Gamma$ is one-to-one. This gives the
needed contradiction.\fdem
\end{demo}

According to \cite{sarkis2} (control of volume), $[\widetilde T]$
is of locally finite mass in $\mathbb{C} \backslash \{0\}$.
Thus according to  Bishop-King-Harvey-Shiffman-Alexander
\cite{bishop,king,harvey2,shiffman,alexander},
$[\widetilde T]$ extends to a solution
(that we will still note $[\widetilde T]$) to the
boundary problem for $[\Gamma_{g,f}]$ in
$\pi^{-1}(\mathbb{C}\backslash \{0\})$.\fdem
\end{demo}

\section{Related problems}

We do not know if theorem $2$ is still valid
if we assume less regularity for $M$ or for the CR functions.
For example,  in the case that $M$ is Lipschitz
and $f$ is in the Sobolev space $W^{-1/2}(M)$,
a counter example
is given by Henkin and Iordan in \cite{iordan}.
Nevertheless, by analogy with the extension result
they obtain, one might expect:
\begin{problem}[Henkin]
Let $\Omega \subset P_n(\mathbb{C})$, $(n \geq 2)$, 
be a domain with lipschitzian boundary $\partial \Omega$
which admits a non constant holomorphic function.
Let $f$ be a CR function which is in the sobolev space 
$W^{1/2}(\partial \Omega)$.
Does $f$ admits a holomorphic
extension in $\mathcal{O}(\Omega) \cap W^1(\Omega)$ ? 
\end{problem}

As we have seen, in the case $M$ is of class $C^2$, the
main difficulty is the possible existence of laminated
compact subset $K$ of $M$. Thus the following problems
become natural:

\begin{problem}
Let $M$ be a connected $C^2$ hypersurface of 
$P_n(\mathbb{C})$ ($n \geq 2$) which divides $P_n(\mathbb{C})$ in two 
connected open sets $\Omega_1$ and $\Omega_2$. Then:
\begin{enumerate}
\item Does there exists $i \in \{1,2\}$ such that smooth CR maps 
$f: M \rightarrow P_1(\mathbb{C})$ extends meromorphically to $\Omega_i$ ?
\item Is $M$ globally minimal ?
\item Does there always exists a non constant $C^1$ CR function on $M$ ?
\end{enumerate}
\end{problem}

Let $U \subset P_n(\mathbb{C})$ be an open set. 
If $U$ contains a laminated compact set $K$ then
holomorphic functions on $U$ have to be constant and meromorphic
functions have to be rational. As we have proved,
continuous CR functions on $K$ are also constant. Thus
one could expect:
\begin{problem}
Let $K \subset P_2(\mathbb{C})$ be a compact set
$C^2$-laminated by Riemann surfaces.
Let $f$ be a $C^2$ CR map from $K$ to $P_1(\mathbb{C})$
(i.e. for any analytic disc $\Delta \subset K$, $f|_\Delta$
is a holomorphic map). Does there exists a rational map $Q: P_2(\mathbb{C}) 
\rightarrow P_1(\mathbb{C})$ such that $Q|_K=f$ ?
\end{problem}

Of course, in the case it is known that there 
exist no non-trivial laminated compact subset of 
$P_2(\mathbb{C})$, the problems $2$ and $3$ 
would be obvious.

\noindent
e-mail: {\tt sarkis@math.jussieu.fr}\par\noindent
\par
\medskip
\end{document}